\documentclass[twocolumn,10pt]{asme2e}

\usepackage{amssymb,epsfig,amsmath,epstopdf}
\usepackage{dsfont}   
\numberwithin{figure}{section} 
\usepackage{graphicx,subfigure}
\newfont{\smalll}{cmr8}
\def\IR{\mathbb{R}}
\def\IS{\mathbb{S}}

\def\zero{\mathbf{0}}
\def\one{\mathbf{1}}

\def\diag{\mathrm{diag}}
\def\eye{\mathrm{I}}

\def\IS{\hbox{I\hskip-.1em S}}
\def\IC{\hbox{C\hskip-
.5em\raise.5ex\hbox{$\scriptscriptstyle\mid$}}\ }
\def\Ic{\hbox{\smalll C\hskip-
.5em\raise.3ex\hbox{$\scriptscriptstyle\mid$}}\ }

\def\T={\buildrel {\scriptscriptstyle\triangle} \over =}

\def\sqr#1#2{{\vcenter{\vbox{\hrule height.#2pt\hbox{\vrule
width.#2pt height#1pt \kern#1pt\vrule width.#2pt}\hrule
height.#2pt}}}}

\def\diag{\mathop{\rm diag}}

\def\block-diag{\mathop{\rm block{\scriptstyle -}diag}}

\def\pmbb#1{\setbox0=\hbox{#1}\raise 0.5ex\box0}

\def\norm#1{\|#1\|}

\newcommand{\bequ}{\begin{eqnarray}}
\newcommand{\eequ}{\end{eqnarray}}
\newcommand{\mT}{^\mathrm{T}}

\newcommand{\rom}{\mathrm}

\newcommand {\teq}      {\triangleq}

\newcommand {\beq}      {\begin{equation}}
\newcommand {\eeq}      {\end{equation}}

\def\IR{{\mathbb R}}
\def\IC{{\mathbb C}}

\def\IS{{\mathbb S}}

\renewcommand\theenumi{\roman{enumi}}
\renewcommand\theenumii{\alph{enumii}}
\renewcommand\theenumiii{\arabic{enumiii}}
\renewcommand\theenumiv{\Alph{enumiv}}

\def\sc{s_{\rm c}}

\def\sc{s_{\rm c}}

\DeclareMathAlphabet{\mathcal}{OMS}{cmsy}{m}{n} 

\title{CONTROL OF MULTIAGENT FORMATIONS\vspace{-0.5cm}}
\author{Dung Tran
    \affiliation{
	Advanced Systems Research Laboratory\\
	Dept. of Mechanical and Aerospace Engineering\\
	Missouri Univ. of Science and Technology\\
	Rolla, Missouri 65409, USA\\
	dmtt24@mst.edu
    }	
}

\author{Tansel Yucelen
    \affiliation{Advanced Systems Research Laboratory\\
    	Dept. of Mechanical and Aerospace Engineering\\
    	Missouri Univ. of Science and Technology\\
    	Rolla, Missouri 65409, USA\\
 	yucelen@mst.edu
    }
}

\begin{document} 
\maketitle

%%%%%%%%%%%%%%%%%%%%%%%%%%%%%%%%%%%%%%%%%%%%%%%%%%%%%%%%%%%%%%%%%%%%%%
%%%%%%%%%%%%%%%%%%%%%%%%%%%%%%%%%%%%%%%%%%%%%%%%%%%%%%%%%%%%%%%%%%%%%%
%%%%%%%%%%%%%%%%%%%%%%%%%%%%%%%%%%%%%%%%%%%%%%%%%%%%%%%%%%%%%%%%%%%%%%
%%%%%%%%%%%%%%%%%%%%%%%%%%%%%%%%%%%%%%%%%%%%%%%%%%%%%%%%%%%%%%%%%%%%%%
%%%%%%%%%%%%%%%%%%%%%%%%%%%%%%%%%%%%%%%%%%%%%%%%%%%%%%%%%%%%%%%%%%%%%%

\begin{abstract}
\noindent {\textit{This paper makes the first attempt to show how information exchange rules represented by a network having multiple layers (multiplex information networks) can be designed for enabling spatially evolving multiagent formations.
Toward this goal, we consider the invariant formation problem and introduce a distributed control architecture that allows capable agents to spatially alter the resulting formation's density and orientation without requiring global information exchange ability. 
%Efficacy of the proposed approach is illustrated on numerical examples.
}}

\baselineskip 16pt

\end{abstract}

\vspace{-0.25cm}

%%%%%%%%%%%%%%%%%%%%%%%%%%%%%%%%%%%%%%%%%%%%%%%%%%%%%%%%%%%%%%%%%%%%%%
%%%%%%%%%%%%%%%%%%%%%%%%%%%%%%%%%%%%%%%%%%%%%%%%%%%%%%%%%%%%%%%%%%%%%%
%%%%%%%%%%%%%%%%%%%%%%%%%%%%%%%%%%%%%%%%%%%%%%%%%%%%%%%%%%%%%%%%%%%%%%
%%%%%%%%%%%%%%%%%%%%%%%%%%%%%%%%%%%%%%%%%%%%%%%%%%%%%%%%%%%%%%%%%%%%%%
%%%%%%%%%%%%%%%%%%%%%%%%%%%%%%%%%%%%%%%%%%%%%%%%%%%%%%%%%%%%%%%%%%%%%%

\section{INTRODUCTION}

\noindent Current distributed control methods (see, for example, \cite{Ren08:01,Mesbahi10,Oh10} and references therein) have a lack of information exchange infrastructure to enable spatially evolving multiagent formations. 
Because, they are designed based on information exchange rules represented by a network having single layer, which lead to multiagent formations having fixed, non-evolving spatial properties. 
For situations where capable agents have to control the resulting formation using these methods, they can only do so if such vehicles have global information exchange ability. 

The contribution of this paper is to make the first attempt to show how information exchange rules represented by a network having multiple layers (multiplex information networks) can be designed for enabling spatially evolving multiagent formations.
Toward this goal, we consider the invariant formation problem [Section 6.3.1, \citen{Mesbahi10}] and introduce a distributed control architecture that allows capable agents to spatially alter the resulting formation's density and orientation without requiring global information exchange ability. 
%Efficacy of the proposed approach is illustrated on numerical examples.

The notation used in this paper is fairly standard. 
Specifically, $\IR$ denotes the set of real numbers, 
$\IR^n$ denotes the set of $n \times 1$ real column vectors, 
$\IR^{n \times m}$ denotes the set of $n \times m$ real matrices, 
$\IR_{+}$ (resp. $\overline{\IR}_{+}$) denotes the set of positive (resp. non-negative-definite) real numbers, 
$\IR^{n \times n}_{+}$ (resp. $\overline{\IR}^{n \times n}_{+}$) denotes the set of $n \times n$ positive-definite (resp. non-negative-definite) real matrices, 
$\IS^{n \times n}_{+}$ (resp. $\overline{\IS}^{n \times n}_{+}$) denotes the set of $n \times n$ symmetric positive-definite (resp. symmetric nonnegative-definite) real matrices, 
$\zero_n$ denotes the $n \times 1$ vector of all zeros,
$\one_n$ denotes the $n \times 1$ vector of all ones,
$\zero_{n \times n}$ denotes the $n \times n$ zero matrix,
and $\textbf{I}_{n}$ denotes the $n \times n$ identity matrix. 
In addition, we write $(\cdot)\mT$ for transpose, 
$\lambda_{\rom{min}}(A)$ (resp. $\lambda_{\rom{max}}(A))$ for the minimum and respectively maximum eigenvalue of the Hermitian matrix A, 
$\lambda_{i}(A)$ for the \textit{i}-th eigenvalue of A (A is symmetric and the eigenvalues are ordered from least to greatest value), 
$\rom{det}(A)$ for the determinant of $A$,
$\diag(a)$ for the diagonal matrix with the vector $a$ on its diagonal, 
$[x]_{i}$ for the entry of the vector $x$ on the \textit{i}-th row, 
and $[A]_{ij}$  for the entry of the of the matrix A on the \textit{i}-th row and \textit{j}-th column.

\vspace{-0.25cm}

%%%%%%%%%%%%%%%%%%%%%%%%%%%%%%%%%%%%%%%%%%%%%%%%%%%%%%%%%%%%%%%%%%%%%%
%%%%%%%%%%%%%%%%%%%%%%%%%%%%%%%%%%%%%%%%%%%%%%%%%%%%%%%%%%%%%%%%%%%%%%
%%%%%%%%%%%%%%%%%%%%%%%%%%%%%%%%%%%%%%%%%%%%%%%%%%%%%%%%%%%%%%%%%%%%%%
%%%%%%%%%%%%%%%%%%%%%%%%%%%%%%%%%%%%%%%%%%%%%%%%%%%%%%%%%%%%%%%%%%%%%%
%%%%%%%%%%%%%%%%%%%%%%%%%%%%%%%%%%%%%%%%%%%%%%%%%%%%%%%%%%%%%%%%%%%%%%

\section{PRELIMINARIES}

\noindent In this section, we recall some of the basic notions from graph theory, which is followed by the general setup of the consensus and formation problems for multiagent systems. 
We refer to \cite{Mesbahi10,godsil2001algebraic} for details about graph theory and multiagent systems.
Specifically, graphs are broadly adopted to encode interactions in networked systems.
An {undirected} graph $\mathcal{G}$ is defined by a set $\mathcal{V}_\mathcal{G}=\{1,\ldots,n\}$ of {nodes}
and a set $\mathcal{E}_\mathcal{G} \subset \mathcal{V}_\mathcal{G} \times \mathcal{V}_\mathcal{G}$ of {edges}.
If $(i,j) \in \mathcal{E}_\mathcal{G}$, then the nodes $i$ and $j$ are {neighbors} and the neighboring relation is indicated with $i \sim j$.
The {degree} of a node is given by the number of its neighbors.
Letting $d_i$ be the degree of node $i$, then the {degree} matrix of a graph $\mathcal{G}$, $\mathcal{D}(\mathcal{G}) \in \IR^{n \times n}$, is given by
\bequ
        \mathcal{D}(\mathcal{G}) \triangleq \diag(d), \quad d=[d_1,\ldots,d_n]\mT. \label{DegMat}
\eequ
A {path} $i_0 i_1 \ldots i_L$ is a finite sequence of nodes such that $i_{k-1} \sim i_k$, $k=1, \ldots, L$,
and a graph $\mathcal{G}$ is {connected} if there is a path between any pair of distinct nodes.
The {adjacency} matrix of a graph $\mathcal{G}$, $\mathcal{A}(\mathcal{G}) \in \IR^{n \times n}$, is given by
\begin{eqnarray}
   [\mathcal{A}(\mathcal{G})]_{ij}
   \teq
   \left\{ \begin{array}{cl}
      1, & \mbox{ if $(i,j)\in\mathcal{E}_\mathcal{G}$},
      \\
      0, & \mbox{otherwise}.
   \end{array} \right.
   \label{AdjMat}
\end{eqnarray}
The {Laplacian} matrix of a graph, $\mathcal{L}(\mathcal{G}) \in \overline{\IS}_+^{\hspace{0.1em} n \times n}$, playing a central role in many graph theoretic treatments of multiagent systems is given by
\bequ
    \mathcal{L}(\mathcal{G}) \triangleq \mathcal{D}(\mathcal{G}) - \mathcal{A}(\mathcal{G}), \label{Laplacian}
\eequ
where the spectrum of the Laplacian for an undirected and connected graph $\mathcal{G}$ can be ordered as $ 0 = \lambda_1(\mathcal{L}(\mathcal{G}))<\lambda_2(\mathcal{L}(\mathcal{G}))\le \cdots \le \lambda_n(\mathcal{L}(\mathcal{G}))$ 
with $\one_n$ as the eigenvector corresponding to the zero eigenvalue $\lambda_1(\mathcal{L}(\mathcal{G}))$
and $\mathcal{L}(\mathcal{G}) \one_n = \zero_n$ and $\rom{e}^{\mathcal{L}(\mathcal{G})} \one_n = \one_n$ hold. 
Throughout this paper, we assume that the graph $\mathcal{G}$ is undirected and connected.

We can model a given multiagent system by a graph $\mathcal{G}$ where nodes and edges represent agents and interagent information exchange links,
respectively. 
Let $x_i(t) \in \IR^m$ denote the state of node $i$, whose dynamics is described by the single integrator
\bequ
        \dot{x}_i(t) = u_i(t), \quad x_i(0)=x_{i0}, \quad i=1,\cdots,n, \label{eqn:01}
\eequ
with $u_i(t) \in \IR^m$ being the control input of node $i$.
Allowing agent $i$ to have access to the relative state information with respect to its neighbors, the solution of the consensus problem can be achieved by applying
\bequ
            u_i(t)=-\sum_{i \sim j} \bigl(x_i(t)-x_j(t) \bigl), \label{eqn:02}
\eequ
to the single integrator dynamics given by (\ref{eqn:01}), where (\ref{eqn:01}) in conjunction with (\ref{eqn:02}) can be represented as the Laplacian dynamics of the form
\bequ
        \dot{x}(t)=-\mathcal{L}(\mathcal{G})\otimes \eye_m \hspace{0.2em} x(t), \quad x(0)=x_0, \label{eqn:03}
\eequ
where $x(t)=[x_1\mT(t),\cdots,x_n\mT(t)]\mT$ denotes the aggregated state vector of the multiagent system. 
Since the graph $\mathcal{G}$ is assumed to be undirected and connected, it follows from (\ref{eqn:03}) that
\bequ
\lim_{t\rightarrow\infty} [x_i(t)]_j&=&\frac{[x_1(0)]_j+\cdots[x_n(0)]_j}{n},
\eequ
holds for $i=1,\ldots,n$ and $j=1,\ldots,m$. 
Throughout this paper, we assume without loss of generality that $m=2$, which implies that the multiagent system evolves in a planar space.

For our take on the formation problem, define $\tau_i(t)\in\IR^2$ as the displacement of $x_i(t)\in\IR^2$ from the target location $\xi_i\in\IR^2$.
Then, using the state transformation given by
\bequ
        \tau_i(t)=x_i(t)-\xi_i, \quad i=1,\ldots,n, \label{eqn:05}
\eequ
the solution of the {invariant} formation problem [Section 6.3.1, \citen{Mesbahi10}] follows from (\ref{eqn:03}) with $m=2$ as
\bequ
        \dot{x}(t)=-\mathcal{L}(\mathcal{G})\otimes \eye_2 \hspace{0.2em}x(t)+\mathcal{L}(\mathcal{G})\otimes \eye_2 \hspace{0.2em}\xi, \quad x(0)=x_0, \label{eqn:06}
\eequ
where $\xi=[\xi_1, \cdots, \xi_n]\mT$.
Note that (\ref{eqn:06}) can equivalently be written as
 \bequ
        \dot{x}_i(t) &=& -\sum_{i \sim j} \bigl(x_i(t)-x_j(t) \bigl)+\sum_{i \sim j}\bigl(\xi_i-\xi_j \bigl), \ x_i(0)=x_{i0}. \ \ \ \ \ \label{eqn:07}
\eequ
In order to present the main results of this paper, we consider this particular formation problem and introduce a distributed control architecture to allow capable agents to spatially alter the resulting formation's density and orientation without requiring global information exchange ability. 

%%%%%%%%%%%%%%%%%%%%%%%%%%%%%%%%%%%%%%%%%%%%%%%%%%%%%%%%%%%%%%%%%%%%%%
%%%%%%%%%%%%%%%%%%%%%%%%%%%%%%%%%%%%%%%%%%%%%%%%%%%%%%%%%%%%%%%%%%%%%%
%%%%%%%%%%%%%%%%%%%%%%%%%%%%%%%%%%%%%%%%%%%%%%%%%%%%%%%%%%%%%%%%%%%%%%
%%%%%%%%%%%%%%%%%%%%%%%%%%%%%%%%%%%%%%%%%%%%%%%%%%%%%%%%%%%%%%%%%%%%%%
%%%%%%%%%%%%%%%%%%%%%%%%%%%%%%%%%%%%%%%%%%%%%%%%%%%%%%%%%%%%%%%%%%%%%%
\vspace{-0.45cm}

\section{\uppercase{Multiplex Information Networks for Invariant Formation Problem}}

\noindent We first introduce a multiplex information networks-based distributed control architecture for controlling density of multiagent formations (Section 3.1) and then generalize our results for controlling both density and orientation of multiagent formations (Section 3.2).

%%%%%%%%%%%%%%%%%%%%%%%%%%%%%%%%%%%%%%%%%%%%%%%%%%%%%%%%%%%%%%%%%%%%%%
%%%%%%%%%%%%%%%%%%%%%%%%%%%%%%%%%%%%%%%%%%%%%%%%%%%%%%%%%%%%%%%%%%%%%%
%%%%%%%%%%%%%%%%%%%%%%%%%%%%%%%%%%%%%%%%%%%%%%%%%%%%%%%%%%%%%%%%%%%%%%
%%%%%%%%%%%%%%%%%%%%%%%%%%%%%%%%%%%%%%%%%%%%%%%%%%%%%%%%%%%%%%%%%%%%%%
%%%%%%%%%%%%%%%%%%%%%%%%%%%%%%%%%%%%%%%%%%%%%%%%%%%%%%%%%%%%%%%%%%%%%%
\vspace{-0.25cm}
\subsection{Formation Density Control}

\noindent Consider a system of $n$ agents exchanging information among each other using their local measurements, according to an undirected and connected graph $\mathcal{G}$. 
In order to control density of the invariant formation problem introduced in the previous section, we propose the distributed controller having two layers given by
\bequ
 	\dot{x}_i(t)&=&-\sum_{i \sim j}\bigl(x_i(t)-x_j(t)\bigl)+\sum_{i \sim j}\bigl( \gamma_i(t)\xi_i -\gamma_j(t)\xi_j \bigl) \nonumber\\
	                 &&-\xi_i \sum_{i \sim j} \bigl(\gamma_i(t)-\gamma_j(t)\bigl)-k_i \xi_i \bigl(\gamma_i(t)-\gamma\bigl), \ \ x_i(0)=x_{i0}, \label{TY:01}\ \ \ \ \ \ \\
	\dot{\gamma}_i(t)&=&-\sum_{i \sim j}\bigl(\gamma_i(t)-\gamma_j(t)\bigl)-k_i\bigl(\gamma_i(t)-\gamma), \ \ \gamma_i(0)=\gamma_{i0}, \label{TY:02}
\eequ
where $x_i(t)\in\IR^2$ denotes the state of the first layer of node $i$ that corresponds to the actual state of node $i$, 
$\xi_i\in\IR^2$ denotes the formation shape of node $i$, 
$\gamma_i(t)\in\IR$ denotes the state of the second layer of node $i$ that is introduced to distribute the formation density parameter $\gamma\in\IR$ through local information exchange, 
and $k_i=1$ for capable (leader) agents and otherwise $k_i=0$ (we implicitly assume that there exists at least one capable agent in the network). 
Note that the formation density parameter $\gamma$ is only available to capable agents as such they have a capability to alter the density of the resulting formation (i.e., scale the formation) through peer-to-peer communications. 

Considering the networked multiagent system given by (\ref{TY:01}) and (\ref{TY:02}), where agents exchange information using local measurements and with $\mathcal{G}$ defining an undirected and connected graph topology, it can be shown that 
\bequ
	\lim_{t\rightarrow\infty} \bigl(x_i(t)-x_j(t) \bigr) &=& \gamma \bigl(\xi_i - \xi_j\bigl), \label{TY:03}
\eequ
holds for all $i=1,\ldots,n$.
This shows that the proposed algorithm given by (\ref{TY:01}) and (\ref{TY:02}) allows the density of the multiagent formation to be controlled by the formation density parameter $\gamma$, which is only available to capable agents. 

Note that the result given by (\ref{TY:03}) can be generalized to the case where the formation density parameter is a bounded function of time with a bounded time rate of change, i.e., $\norm{\gamma(t)}_2\le\gamma^*$ and $\norm{\dot{\gamma}(t)}_2\le\dot{\gamma}^*$. In this case, however, (\ref{TY:03}) needs to be replaced with 
\bequ
\lim_{t\rightarrow\infty} \bigl(x_i(t)-x_j(t)\bigr) &=& \gamma_i(t)\xi_i-\gamma_j(t)\xi_j, \label{TY:18}
\eequ 
where $\gamma_i(t)$ and $\gamma_j(t)$ converge to a neighborhood of $\gamma(t)$. 
In order to make this neighborhood close to $\gamma(t)$, it can be assumed that $\dot{\gamma}^*$ is small.  If this is not a valid assumption depending on an application of interest, then one can consider the distributed controller given by
\bequ
 	\dot{x}_i(t)&=&-\sum_{i \sim j}\bigl(x_i(t)-x_j(t)\bigl)+\sum_{i \sim j}\bigl( \gamma_i(t)\xi_i -\gamma_j(t)\xi_j \bigl) \nonumber\\
	                 &&-\alpha\xi_i \sum_{i \sim j} \bigl(\gamma_i(t)-\gamma_j(t)\bigl)-\alpha k_i \xi_i \bigl(\gamma_i(t)-\gamma\bigl), \label{TY:21}\ \ \ \ \ \ \\
	\dot{\gamma}_i(t)&=&-\alpha\sum_{i \sim j}\bigl(\gamma_i(t)-\gamma_j(t)\bigl)-\alpha k_i\bigl(\gamma_i(t)-\gamma), \label{TY:22}
\eequ
and increase $\alpha\in\IR_{+}$ in order to drive $\gamma_i(t)$ and $\gamma_j(t)$ to a close neighborhood of $\gamma(t)$.

%%%%%%%%%%%%%%%%%%%%%%%%%%%%%%%%%%%%%%%%%%%%%%%%%%%%%%%%%%%%%%%%%%%%%%
%%%%%%%%%%%%%%%%%%%%%%%%%%%%%%%%%%%%%%%%%%%%%%%%%%%%%%%%%%%%%%%%%%%%%%
%%%%%%%%%%%%%%%%%%%%%%%%%%%%%%%%%%%%%%%%%%%%%%%%%%%%%%%%%%%%%%%%%%%%%%
%%%%%%%%%%%%%%%%%%%%%%%%%%%%%%%%%%%%%%%%%%%%%%%%%%%%%%%%%%%%%%%%%%%%%%
%%%%%%%%%%%%%%%%%%%%%%%%%%%%%%%%%%%%%%%%%%%%%%%%%%%%%%%%%%%%%%%%%%%%%%
\vspace{-0.25cm}
\subsection{Formation Density and Orientation Control}

\noindent As in Section 3.1, consider a system of $n$ agents exchanging information among each other using their local measurements, according to an undirected and connected graph $\mathcal{G}$. In order to control both density and orientation of the invariant formation problem introduced in the previous section, we propose the distributed controller having three layers given by
\bequ
\dot{x}_i(t) &=& -\sum_{i \sim j}\bigl(x_i(t) -x_j(t)\bigr) \nonumber \\
&&+\sum_{i \sim j}\Bigl(\gamma_i(t)R(\theta_i(t))\xi_i- \gamma_j(t)R(\theta_j(t))\xi_j\Bigl) \nonumber \\
&& -\Bigl(\sum_{i \sim j}\bigl(\gamma_i(t)- \gamma_j(t)\bigr) + k_i\bigl(\gamma_i(t) - \gamma\bigr)\Bigr)R(\theta_i(t))\xi_i \nonumber \\
&&-\gamma_i(t)\Bigl(\sum_{i \sim j}\bigl(\theta_i(t) - \theta_j(t)\bigr)+ k_i\bigl(\theta_i(t) - \theta\bigr)\Bigr)Q(\theta_i(t))\xi_i, \nonumber \\
&& \hspace{5.3cm} x_{i}(0) = x_{i0}, \label{DT:01} \\
\dot{\gamma}_i(t) &=& -\sum_{i \sim j}\bigl(\gamma_i(t)- \gamma_j(t)\bigr)- k_i\bigl(\gamma_i(t) - \gamma\bigr),\quad \gamma_i(0) = \gamma_{i0},  \label{DT:02}\\
\dot{\theta}_i(t) &=& -\sum_{i \sim j}\bigl(\theta_i(t) - \theta_j(t)\bigr)- k_i\bigl(\theta_i(t) - \theta\bigr),\hspace{0.2cm} \theta_i(0) = \theta_{i0}, \label{DT:03} 
\eequ
where $x_i(t)\in\IR^2$ denotes the state of the first layer of node $i$ that corresponds to the actual state of node $i$, 
$\xi_i\in\IR^2$ denotes the formation shape of node $i$, 
$\gamma_i(t)\in\IR$ denotes the state of the second layer of node $i$ that is introduced to distribute the formation density parameter $\gamma\in\IR$ through local information exchange, $\theta_i(t) \in \IR$ denotes the state of the third layer of node $i$ that is introduced to distribute the formation orientation parameter $\theta \in \IR$ through local information exchange, and $k_i=1$ for capable (leader) agents and otherwise $k_i=0$. 
In (\ref{DT:01}), $R(\theta_i(t))$ denotes the rotation matrix of agent $i$
\bequ
R(\theta_i(t)) \triangleq \begin{bmatrix} 
	\cos{\theta_i(t)} & -\sin{\theta_i(t)} \\
	\sin{\theta_i(t)} & \cos{\theta_i(t)} 
\end{bmatrix},  \label{RotMat}
\eequ
and 
\vspace{0cm}
\bequ
Q(\theta_i(t)) \triangleq \frac{dR(\theta_i(t))}{d\theta_i(t)}.
\eequ 
%\bequ
%Q(\theta_i(t)) \triangleq \frac{dR(\theta_i(t))}{d\theta_i(t)} = \begin{bmatrix} 
%	-\sin{\theta_i(t)} & -\cos{\theta_i(t)} \\
%	\cos{\theta_i(t)} & -\sin{\theta_i(t)} 
%\end{bmatrix}, \quad \label{DRotMatx}
%\eequ
As in the previous section, note that the formation density parameter $\gamma$ and the orientation parameter $\theta$  are only available to capable agents as such they have the capability to alter the density and orientation of the resulting formation (i.e., scale and rotate the formation). 

Considering the networked multiagent system given by (\ref{DT:01}), (\ref{DT:02}), and (\ref{DT:03}), where agents exchange information using local measurements and with $\mathcal{G}$ defining an undirected and connected graph topology, it can be shown that 
\bequ
\lim_{t\rightarrow\infty} \bigl(x_i(t)-x_j(t)\bigr) &=& \gamma R(\theta)\bigl(\xi_i - \xi_j\bigl), \label{DT:04}
\eequ
holds for all $i=1,\ldots,n$. 
This result shows that the proposed algorithm given by (\ref{DT:01}), (\ref{DT:02}) and (\ref{DT:03}) allows the density and orientation of the multiagent formation to be controlled by the formation density parameter $\gamma$ and orientation parameter $\theta$, which are only available to capable agents.

Note that the result given by ({\ref{DT:04}) can be generalized to the case where the formation density and orientation parameters are bounded functions of time with bounded time rates of change, i.e., $\norm{\gamma(t)}_2\le\gamma^*$, $\norm{\theta(t)}_2\le\theta^*$, $\norm{\dot{\gamma}(t)}_2\le\dot{\gamma}^*$, and $\norm{\dot{\theta}(t)}_2\le\dot{\theta}^*$. In this case, however, (\ref{DT:04}) needs to be replaced with 
\bequ
\lim_{t\rightarrow\infty} \bigl(x_i(t)-x_j(t)\bigr) &=& \gamma_i(t)R(\theta_i(t))\xi_i - \gamma_j(t)R(\theta_j(t))\xi_j, \ \ \ \ \label{DT:16}
\eequ
where $\gamma_i(t)$ and $\gamma_j(t)$ converge to a neighborhood of $\gamma(t)$ and $\theta_i(t)$ and $\theta_j(t)$  converge to a neighborhood of $\theta(t)$. In order to make this neighborhood close to $\gamma(t)$ and  $\theta(t)$ , it can be assumed that $\dot{\gamma}^*$ and $\dot{\theta}^*$ are small.  If this is not a valid assumption, then one can consider the distributed controller given by
\bequ
\dot{x}_i(t) &=& -\sum_{i \sim j}\bigl(x_i(t) -x_j(t)\bigr) \nonumber \\
&&+\sum_{i \sim j}\Bigl(\gamma_i(t)R(\theta_i(t))\xi_i- \gamma_j(t)R(\theta_j(t))\xi_j\Bigl) \nonumber \\
&& -\alpha\Bigl(\sum_{i \sim j}\bigl(\gamma_i(t)- \gamma_j(t)\bigr) + k_i\bigl(\gamma_i(t) - \gamma\bigr)\Bigr)R(\theta_i(t))\xi_i \nonumber \\
&&-\alpha\gamma_i(t)\Bigl(\sum_{i \sim j}\bigl(\theta_i(t) - \theta_j(t)\bigr)+ k_i\bigl(\theta_i(t) - \theta\bigr)\Bigr)Q(\theta_i(t))\xi_i, \nonumber\\ \label{DT:17}\\
\dot{\gamma}_i(t) &=& -\alpha\sum_{i \sim j}\bigl(\gamma_i(t)- \gamma_j(t)\bigr)- \alpha k_i\bigl(\gamma_i(t) - \gamma\bigr),  \label{DT:18}\\
\dot{\theta}_i(t) &=& -\alpha\sum_{i \sim j}\bigl(\theta_i(t) - \theta_j(t)\bigr)- \alpha k_i\bigl(\theta_i(t) - \theta\bigr), \label{DT:19} 
\eequ
and increase $\alpha\in\IR_{+}$ in order to drive $\gamma_i(t)$ and $\gamma_j(t)$ to a close neighborhood of $\gamma(t)$ as well as $\theta_i(t)$ and $\theta_j(t)$ to a close neighborhood of $\theta(t)$. 

%%%%%%%%%%%%%%%%%%%%%%%%%%%%%%%%%%%%%%%%%%%%%%%%%%%%%%%%%%%%%%%%%%%%%%
%%%%%%%%%%%%%%%%%%%%%%%%%%%%%%%%%%%%%%%%%%%%%%%%%%%%%%%%%%%%%%%%%%%%%%
%%%%%%%%%%%%%%%%%%%%%%%%%%%%%%%%%%%%%%%%%%%%%%%%%%%%%%%%%%%%%%%%%%%%%%
%%%%%%%%%%%%%%%%%%%%%%%%%%%%%%%%%%%%%%%%%%%%%%%%%%%%%%%%%%%%%%%%%%%%%%
%%%%%%%%%%%%%%%%%%%%%%%%%%%%%%%%%%%%%%%%%%%%%%%%%%%%%%%%%%%%%%%%%%%%%%
\vspace{-0.25cm}
\section{\uppercase{Conclusion}}
\noindent To contribute to the previous studies in formation control of multiagent systems, we considered the invariant formation problem and presented a multiplex information networks-based distributed control architecture.
The proposed methodology allows capable agents to spatially alter the resulting formation's density and orientation without requiring global information exchange ability. 
%Numerical examples illustrated the efficacy of the proposed approach. 

\vspace{-0.35cm}
%%%%%%%%%%%%%%%%%%%%%%%%%%%%%%%%%%%%%%%%%%%%%%%%%%%%%%%%%%%%%%%%%%%%%%
%%%%%%%%%%%%%%%%%%%%%%%%%%%%%%%%%%%%%%%%%%%%%%%%%%%%%%%%%%%%%%%%%%%%%%
%%%%%%%%%%%%%%%%%%%%%%%%%%%%%%%%%%%%%%%%%%%%%%%%%%%%%%%%%%%%%%%%%%%%%%
%%%%%%%%%%%%%%%%%%%%%%%%%%%%%%%%%%%%%%%%%%%%%%%%%%%%%%%%%%%%%%%%%%%%%%
%%%%%%%%%%%%%%%%%%%%%%%%%%%%%%%%%%%%%%%%%%%%%%%%%%%%%%%%%%%%%%%%%%%%%%

\bibliographystyle{asmems4}
\bibliography{MultiplexBib}
\end{document}